\documentclass[11pt]{article}
\usepackage{amssymb,amsmath}
\usepackage{amsmath}
\usepackage{amsfonts}

\newtheorem{precor}{{\bf Corollary}}

\newtheorem{precon}{{\bf Conjecture}}

\newtheorem{predefin}{{\bf Definition}}

\newtheorem{preexm}{{\bf Example}}

\newtheorem{preappl}{{\bf Application}}

\newtheorem{prelem}{{\bf Lemma}}

\newtheorem{preproof}{{\bf Proof.\ }}

\newenvironment{proof}[1]{\begin{preproof}{\rm
               #1}\hfill{$\blacksquare$}}{\end{preproof}}
\newtheorem{presproof}{{\bf Sketch of Proof.\ }}

\newtheorem{prethm}{{\bf Theorem}}

\newtheorem{prealphthm}{{\bf Theorem}}

\newtheorem{prepro}{{\bf Proposition}}

\newtheorem{preprb}{{\bf Problem}}

\def\conct[#1,#2]{\mbox {${#1} \leftrightarrow {#2}$}}
\def\dconct[#1,#2]{\mbox {${#1} \rightarrow {#2}$}}
\def\deg[#1,#2]{\mbox {$d_{_{#1}}(#2)$}}
\def\mindeg[#1]{\mbox {$\delta_{_{#1}}$}}
\def\maxdeg[#1]{\mbox {$\Delta_{_{#1}}$}}
\def\outdeg[#1,#2]{\mbox {$d_{_{#1}}^{^+}(#2)$}}
\def\minoutdeg[#1]{\mbox {$\delta_{_{#1}}^{^+}$}}
\def\maxoutdeg[#1]{\mbox {$\Delta_{_{#1}}^{^+}$}}
\def\indeg[#1,#2]{\mbox {$d_{_{#1}}^{^-}(#2)$}}
\def\minindeg[#1]{\mbox {$\delta_{_{#1}}^{^-}$}}
\def\maxindeg[#1]{\mbox {$\Delta_{_{#1}}^{^-}$}}

\def\dre[#1,#2,#3]{\mbox {${\cal E}_{_{#3}}(#1,#2)$}}
\def\pdre[#1,#2,#3]{\mbox {${\cal P}_{_{#3}}(#1,#2)$}}
\def\var[#1,#2]{\mbox {${\rm Var}_{_{#1}}(#2)$}}
\def\ls[#1]{\mbox {$\xi^{^{#1}}$}}
\def\hom[#1,#2]{\mbox {${\rm Hom}({#1},{#2})$}}
\def\onvhom[#1,#2]{\mbox {${\rm Hom^{v}}(#1,#2)$}}
\def\onehom[#1,#2]{\mbox {${\rm Hom^{e}}(#1,#2)$}}
\def\core[#1]{\mbox {$#1^{^{\bullet}}$}}
\def\cay[#1,#2]{\mbox {${\rm Cay}({#1},{#2})$}}
\def\cays[#1,#2]{\mbox {${\rm Cay_{s}}({#1},{#2})$}}
\def\dirc[#1]{\mbox {$\stackrel{\rightarrow}{C}_{_{#1}}$}}
\def\cycl[#1]{\mbox {${\bf Z}_{_{#1}}$}}

\def\sdg[#1]{\mbox {$\stackrel{\leftrightarrow}{#1}$}}

\begin{document}
\footnotetext[1]{This paper is partially supported by Shahid
Beheshti University.}
\begin{center}
{\Large \bf A Note on Partial List Coloring}\\
\vspace*{0.5cm}
{\bf Moharram N. Iradmusa$^1$}\\
{\it Department of Mathematical Sciences}\\
{\it Shahid Beheshti University}\\
{\it P.O. Box {\rm 19834}, Tehran, Iran}\\
{\tt iradmusa@yahoo.com}\\
\end{center}
\begin{abstract}
\noindent Let $G$ be a simple graph with $n$ vertices and list
chromatic number $\chi_\ell(G)=\chi_\ell$. Suppose that $0\leq t\leq
\chi_\ell$ and each vertex of $G$ is assigned a list of $t$ colors.
Albertson, Grossman and Haas [1] conjectured that at least
$\frac{tn}{\chi_\ell}$ vertices of $G$ can be colored from these lists.\\
\noindent In this paper we find some new results in partial list
coloring which help us to show that the conjecture is true for at
least half of the numbers of the set
$\{1,2,\ldots,\chi_\ell(G)-1\}$. In addition we introduce a new
related conjecture and finally we present some results about this
conjecture.\\

\begin{itemize}
\item[]{{\footnotesize {\bf Key words:}\ List chromatic number, List coloring, List assignment.}}
\item[]{ {\footnotesize {\bf Subject classification: 05C15} .}}
\end{itemize}
\end{abstract}

\section{Introduction and Preliminaries}

\indent Throughout this paper we only consider finite and simple
graphs. Let $G$ be a simple graph with vertex set $V(G)$, $|V(G)|=n$
and edge set $E(G)$. For each vertex $v\in V(G)$ let
$\mathcal{L}(v)$ be a list of allowed colors assigned to $v$. The
collection of all lists is called a \emph{list assignment} and
denoted by $\mathcal{L}$. In special case we have $t-$\emph{list
assignment} if $|\mathcal{L}(v)|=t$ for all $v\in V(G)$. Also we
call
$R(\mathcal{L})=\bigcup_{v\in V(G)}{\mathcal{L}(v)}$ the \emph{color list} of $\mathcal{L}$.\\
\indent The graph $G$ is called \emph{$\mathcal{L}-$list colorable}
if there is a coloring $c:V(G)\rightarrow R(\mathcal{L})$ such that
$c(v)\neq c(u)$ for all $uv\in E(G)$ and $c(v)\in \mathcal{L}(v)$
for all $v\in V(G)$. Furthermore $G$ is \emph{$k-$choosable} if it
is $\mathcal{L}-$list colorable for every $k-$list assignment
$\mathcal{L}$. The \emph{list chromatic number} or \emph{choice
number}, $\chi_\ell(G)=\chi_\ell$ is the smallest number $k$ such
that $G$ is $k-$choosable. List coloring was
introduced independently by Vising $[6]$ and by Erdos, Rubin and Taylor $[3]$.\\
\indent In addition, let $\lambda_\mathcal{L}(G)$ be the maximum
number of vertices of $G$ which are colorable with respect to the
list assignment $\mathcal{L}$. Define
$\lambda_t(G)=min{\lambda_\mathcal{L}(G)}$, where the minimum is
taken over all $t-$list assignments $\mathcal{L}$. Clearly if $t\geq
\chi_\ell(G)$ then $\lambda_t=n$. Thus it is interesting to know
about $\lambda_t$ when
$t<\chi_\ell$.\\
\textbf{Conjecture 1. (Albertson, Grossman and Haas [1])} Let $G$ be
a graph with $n$ vertices and list chromatic number
$\chi_\ell(G)=\chi_\ell$. Then, for any $0\leq t\leq \chi_\ell$ we have $\lambda_t\geq\frac{tn}{\chi_\ell}$.\\
\noindent We call this conjecture \emph{AGH conjecture}. The
conjecture is clearly correct for $t=0$, $t=1$ and $t=\chi_\ell$.
Albertson et. al. [1] proved the following theorem:\\
\textbf{Theorem 1.} Let $0\leq t\leq \chi_\ell$. Then
$\lambda_{t}\geq (1-(1-\frac{1}{\chi})^t)n$ and this number is
asymptotically best possible.\\
Furthermore, Chappell [2] found the lower bound
$\frac{6}{7}\frac{tn}{\chi_\ell}$ for $\lambda_{t}$ for all $t$
between $0$ and $\chi_\ell$. In addition Haas et. al. [4] showed
that the conjecture is true when $t|\chi_\ell$. For more information see also [5] and [7].\\
In this paper we find some new results which help us to show that
\emph{AGH conjecture} is true for at least half of the numbers of
the set $\{1,2,\ldots,\chi_\ell(G)-1\}$. Finally we introduce a new
conjecture and present some results about this conjecture.
\section{Main Results}
In the first theorem of this section we prove a fundamental
inequality about the
numbers $\lambda_t$ where $0\leq t\leq\chi_\ell$.\\
\textbf{Theorem 2.(Triangle Inequality)} Let $1\leq r,s\leq
\chi_\ell$. Then,
\[\lambda_r+\lambda_s\geq\lambda_{r+s}.\]
\begin{proof}
{Let $\mathcal{L}_r$ and $\mathcal{L}_s$ be $r-$list and $s-$list
assignments of $G$ respectively, such that
$\lambda_{\mathcal{L}_r}=\lambda_r$ and
$\lambda_{\mathcal{L}_s}=\lambda_s$. Also suppose that
$R(\mathcal{L}_r)\cap R(\mathcal{L}_s)=\varnothing$. We define an
$(r+s)-$list assignment $\mathcal{L}_{r+s}$ on $V(G)$ with
\[\mathcal{L}_{r+s}(v)=\mathcal{L}_{r}(v)\cup\mathcal{L}_{s}(v).\]
Then we have $\lambda_{r+s}\leq\lambda_{\mathcal{L}_{r+s}}$. Now
suppose that $c:V(G)\rightarrow R(\mathcal{L}_{r+s})$ is a partial
list coloring of $G$ such that $\lambda_{\mathcal{L}_{r+s}}$
vertices of $G$ are colored. Then we can divide colored vertices of
$G$ to the sets $R=\{v\in V(G)| c(v)\in \mathcal{L}_r(v)\}$ and
$S=\{v\in V(G)| c(v)\in \mathcal{L}_s(v)\}$. This partition results
$|R|\leq \lambda_{\mathcal{L}_r}$ and $|S|\leq
\lambda_{\mathcal{L}_s}$. Finally,
\begin{center}
$\lambda_{r+s}\leq \lambda_{\mathcal{L}_{r+s}}=|R|+|S|\leq
 \lambda_{\mathcal{L}_r}+\lambda_{\mathcal{L}_s}=\lambda_r+\lambda_s$.
\end{center}
}
\end{proof}
\textbf{Corollary 3.} Let $G$ be a graph with list chromatic number
$\chi_\ell(G)=\chi_\ell$.
\begin{enumerate}
\item If $1\leq r_i\leq \chi_\ell(1\leq i\leq k)$, then
$\sum_{i=1}^{k}{\lambda_{r_i}}\geq\lambda_{\sum_{i=1}^{k}{r_i}}$.\\
In particular case, if $r_1=r_2=\ldots=r_k=r$, then we have
$k\lambda_{r}\geq \lambda_{kr}$.
\item If $1\leq r,s\leq \chi_\ell$ and $r|s$, then
$\frac{s}{r}\lambda_{r}\geq \lambda_{s}$ or
$\frac{\lambda_{r}}{r}\geq \frac{\lambda_{s}}{s}$.
\end{enumerate}
The next Corollary was proved in [4]. Here we give another proof.\\
\noindent \textbf{Corollary 4.} If $t|\chi_\ell(G)$, then
$\lambda_t\geq\frac{tn}{\chi_\ell}$.
\begin{proof}
{In part two of Corollary 3 set $s=\chi_\ell$ and use
$\lambda_{\chi_\ell}=|V(G)|$.
 }
\end{proof}
As a result of Theorem 2 we can prove that \emph{AGH conjecture} is
true for at least half of the numbers of the set
$\{1,2,\ldots,\chi_\ell(G)-1\}$.\\
\noindent \textbf{Corollary 5.} Let $1\leq r\leq \chi_\ell=s$. Then
\emph{AGH conjecture} is true for $r$ or $s-r$. In other words, at
least one of the following inequalities is true:\\
\[\lambda_r\geq\frac{rn}{\chi_\ell}\ \ ,\ \ \lambda_{s-r}\geq\frac{(s-r)n}{\chi_\ell}.\]
\begin{proof}
{ Suppose that $\lambda_r<\frac{rn}{\chi_\ell}$ and
$\lambda_{s-r}<\frac{(s-r)n}{\chi_\ell}$. Then
\[\lambda_r+\lambda_{s-r}<\frac{rn}{\chi_\ell}+\frac{(s-r)n}{\chi_\ell}=n.\]
Furthermore from Theorem 2 we have
$\lambda_r+\lambda_{s-r}\geq\lambda_{r+s-r}=\lambda_{s}$.\\
So we conclude $\lambda_{\chi_\ell}<n$ that is a contradiction.}
\end{proof}
Similarly we can prove the following result:\\
\noindent \textbf{Corollary 6.} Let $1\leq r\leq \chi_\ell$ and
$\chi_\ell=kr+r_0$ where $0\leq r_0<r$. Then \emph{AGH conjecture}
is true for $r$ or $r_0$.\\
The next Corollary was proved in [4]. Here we prove it as a result of Corollary 3.\\
\noindent \textbf{Corollary 7.} If $1\leq r\leq \chi_\ell=s$ then
\[\lambda_r\geq\frac{n}{\lceil\frac{s}{r}\rceil}.\]
\begin{proof}
{From Corollary 3 we have
\[\lceil\frac{s}{r}\rceil \lambda_r\geq \lambda_{r\lceil\frac{s}{r}\rceil}.\]
Also from $r\lceil\frac{s}{r}\rceil\geq s$ we have
$\lambda_{r\lceil\frac{s}{r}\rceil}=n$ which completes the proof.}
\end{proof}
In the next corollary we generalize the Corollaries 4 and 6.\\
\noindent \textbf{Corollary 8.} Let $1\leq r\leq s\leq \chi_\ell$.
\begin{enumerate}
\item If $r|s$ and \emph{AGH conjecture} is true for $s$, then it is
true for $r$.
\item Suppose that $s=kr+r_0$ and $0\leq r_0<r$. If \emph{AGH conjecture} is true for
$s$, then it is true for $r$ or $r_0$.
\end{enumerate}
\noindent Here we give a new conjecture that is a generalization of
\emph{AGH conjecture}:\\
\noindent \textbf{Conjecture 2.} Let $G$ be a graph of order $n$
with list chromatic number $\chi_\ell$ and $1\leq r\leq s\leq
\chi_\ell$. Then
\[\frac{\lambda_r}{r}\geq\frac{\lambda_s}{s}.\]
\noindent \textbf{Remark.}
\begin{enumerate}
\item Part two of the Corollary 3
shows that Conjecture 2 is true when $r|s$.
\item If $s=\chi_\ell$, then $\lambda_s=n$. So \emph{AGH
conjecture} is the special case of Conjecture 2.
\item Suppose that $r=1$. So $r|s$ and we have $\lambda_s\leq
s\lambda_1=s\alpha(G)$ where $\alpha(G)$ is the stability number of
$G$.
\end{enumerate}
In the next theorems we show some results about the new conjecture.
We need some definitions which are given before the theorems.\\
\noindent \textbf{Definition 1.} Let $G$ be a graph of order $n$
with list chromatic number $\chi_\ell$. Then $Td(G)$ is the set of
all ordered pairs $(r,s)$, such that $r\leq s$ and
$\frac{\lambda_r}{r}\geq\frac{\lambda_s}{s}$.\\
\noindent \textbf{Theorem 9.} Let $G$ be a graph of order $n$ with
list chromatic number $\chi_\ell$ and $1\leq r\leq s\leq \chi_\ell$.
\begin{enumerate}
\item We have $(r,s)\in Td(G)$ or $(s-r,s)\in Td(G)$.
So Conjecture 2 is true for at least half of the elements of the set
$Td(G)$.
\item If $s=kr+r_0$ and $0\leq
r_0<r$, then $(r,s)\in Td(G)$ or $(r_0,s)\in Td(G)$.
\item We have
\[\lceil\frac{s}{r}\rceil\lambda_r\geq \lambda_s.\]
In other word $\alpha\frac{\lambda_r}{r}\geq\frac{\lambda_s}{s}$,
where $\alpha=\frac{r}{s}\lceil\frac{s}{r}\rceil$.
\item The set $Td(G)$ is an order relation that is defined on the set
$\{1,2,\ldots,\chi_\ell(G)\}$.
\end{enumerate}
\begin{proof}
{We can prove this theorem similar to the Corollaries 5,6 and 7.}
\end{proof}
\noindent \textbf{Definition 2.} Let $\mathcal{L}_1$ and
$\mathcal{L}_2$ be two list assignments of $G$. We say
$\mathcal{L}_1\subseteq \mathcal{L}_2$ if
$\mathcal{L}_1(v)\subseteq\mathcal{L}_2(v)$ for all $v\in V(G)$.\\
\noindent \textbf{Definition 3.} Let $H$ be a subgraph of the graph
$G$, $\mathcal{L}'$ be a list assignment of $H$ and $\mathcal{L}$ be
a list assignment of $G$. We say $\mathcal{L}'$ is a
\emph{restriction of $\mathcal{L}$ to $H$} if
$\mathcal{L}'(v)=\mathcal{L}(v)$ for all $v\in V(H)$ and we show
that by $\mathcal{L}'=\mathcal{L}\mid_H$.\\
The next lemma is an straightforward result of the definition 3.\\
\noindent \textbf{Lemma 10.} Let $H$ be an induced subgraph of the
graph $G$, $\mathcal{L}$ be a list assignment of $G$ and
$\mathcal{L}'=\mathcal{L}\mid_H$. Then $\lambda_{\mathcal{L}}(G)\geq
\lambda_{\mathcal{L}'}(H)$.\\
\noindent \textbf{Theorem 11.} Let $H$
be an induced subgraph of the graph $G$ and $1\leq t\leq
\chi_\ell(G)$. Then $\lambda_{t}(G)\geq \lambda_{t}(H)$.
\begin{proof}
{Suppose that $\mathcal{L}$ is a $t-$list assignment of $G$ such
that $\lambda_{\mathcal{L}}=\lambda_{t}(G)$ and
$\mathcal{L}'=\mathcal{L}\mid_H$. So
$\lambda_{t}(G)=\lambda_{\mathcal{L}}(G)\geq
\lambda_{\mathcal{L}'}(H)\geq \lambda_{t}(H)$. }
\end{proof}
Now we prove the following theorem for the conjecture 2, that is
similar to Theorem 1.\\
\noindent \textbf{Theorem 12.} Let $G$ be a graph of order $n$ with
list chromatic number $\chi_\ell$ and $1\leq r\leq s\leq \chi_\ell$.
Also suppose that $\mathcal{L}_s$ is an $s-$list assignment of $G$
with $\mathcal{L}_s(v)=\{1,2,\ldots,s\}$ for all $v\in V(G)$. Then
\[\lambda_r\geq (1-(1-\frac{1}{s})^r)\lambda_{\mathcal{L}_s}.\]
\begin{proof}
{We can color at most $\lambda_{\mathcal{L}_s}$ vertices of $G$ from
the list color $\{1,2,\ldots,s\}$. Let $H$ be a subgraph of $G$
induced by $\lambda_{\mathcal{L}_s}$ colored vertices. At first we
show that $\chi(H)=s$. Suppose that $\chi(H)\leq s-1$. So there is a
coloring of $G$ by at most $s-1$ colors of $\{1,2,\ldots,s\}$. Thus
we can color at least one vertex of $V(G)-V(H)$ by the remaining
colors such that the coloring remain proper which contradict by the
definition of $H$. So $\chi(H)=s$. Now from Theorem 1 we have
\[\lambda_r(H)\geq
(1-(1-\frac{1}{\chi(H)})^r)|V(H)|=(1-(1-\frac{1}{s})^r)\lambda_{\mathcal{L}_s}.\]
Also from Theorem 11 we have $\lambda_{r}(G)\geq \lambda_{r}(H)$. So
$\lambda_r(G)\geq (1-(1-\frac{1}{s})^r)\lambda_{\mathcal{L}_s}$.}
\end{proof}
\noindent \textbf{Corollary 13.} By assumptions of Theorem 12 and by
use of $\lambda_{\mathcal{L}_s}\geq \lambda_{s}$, we have
\[\lambda_r\geq
(1-(1-\frac{1}{s})^r)\lambda_{s},\]
that is similar to the result of Theorem 1.\\
The last theorem shows the relation between partial list coloring of
a graph and its subgraphs.\\
\noindent \textbf{Theorem 14.} Let $G$ be a graph of order $n$ with
list chromatic number $\chi_\ell$ and $1\leq t\leq \chi_\ell$.\\
(1) Suppose that $H$ is a subgraph of $G$ of maximum order such that
$\chi_\ell(H)=t$. Then, $\lambda_t(G)\geq |H|$.\\
(2) Suppose that
$\lambda_t=\lambda_{\mathcal{L}_t}$ where $\mathcal{L}_t$ is a
$t-$list assignment of $G$. If $H$ is a subgraph of $G$ of order
$\lambda_t$, induced by colored vertices in the partial coloring of
$G$ by $\mathcal{L}_t$, then $\chi_\ell(H)\geq t$.
\begin{proof}
{ (1) Suppose that $\lambda_t(G)\leq |H|-1$ and
$\lambda_t(G)=\lambda_{\mathcal{L}_t}(G)$ where $\mathcal{L}_t$ is a
$t-$list assignment of $G$. Thus in the coloring of $G$ from the
lists of $\mathcal{L}_t$, we can color less than $|H|$ vertices. But
$\chi_\ell(H)=t$. So by restriction of $\mathcal{L}_t$ to $H$, we
can color all vertices of $H$ that is a partial list coloring of $G$
and this contradicts with  $\lambda_{\mathcal{L}_t}(G)<|H|$. So
$\lambda_t(G)\geq |H|$.\\
(2) Suppose that $\chi_\ell(H)=s\leq t-1$. So for any $s-$list
assignment $\mathcal{L}_s$ of $H$, there is a proper coloring of $H$
by the lists of $\mathcal{L}_s$. Also consider the $t-$list
assignment $\mathcal{L}_t$ and select arbitrary vertex $x$ of $G-H$.
Because $x$ is uncolored, so each color of $\mathcal{L}_t(x)$ must
be appeared on at least one vertex of $N(x)$. Now define an $s-$list
assignment $\mathcal{L}'$ such that $\mathcal{L}'\subseteq
\mathcal{L}$ and $c\notin \mathcal{L}'(v)$ for all $v\in V(G)$. So
we can color the vertices of $H$ by lists of $\mathcal{L}'$. Thus by
adding $t-s$ removed colors to each list we can color the vertex $x$
by the color $c$. So we can color $|H|+1$ vertices of $G$ by lists
of $\mathcal{L}_t$ that is a contradiction. So $\chi_\ell(H)\geq
t$.}
\end{proof}
\noindent {\bf Acknowledgment}\\
We would like to thank Professor Hossein Hajiabolhassan for his
useful comments.


\end{document}